\documentclass{article}

\usepackage{amsmath,amssymb,booktabs,xfrac,amsthm,bbm, comment, subfigure, breakcites}

\RequirePackage[OT1]{fontenc}
\RequirePackage[colorlinks,citecolor=blue,urlcolor=blue]{hyperref}

\newtheorem{theorem}{Theorem}

\newtheorem{corollary}{Corollary}
\theoremstyle{remark}
\newtheorem{remark}{Remark}
\theoremstyle{definition}
\newtheorem{definition}{Definition}

\newtheorem*{notation*}{Notation}
\newtheorem*{question*}{Question}
\newtheorem{example}{Example}

\def\y{{\mathbf{y}}}
\def\Y{{\mathbf{Y}}}
\def\Var{{\text{Var}}}
\def\trace{{\text{trace}}}
\def\E{{\text{E}}}

\def\Vstar{V^\star}
\def\Dstar{D}
\def\Vstar{V}

\def\Vtrue{V_{true}}

\def\D0{D_0}
\def\V0{V_0}

\title{Bernstein - von Mises theorem  and misspecified models: a review}

\author{Natalia Bochkina,\\ University of Edinburgh and Maxwell Institute, UK }

\begin{document}

\maketitle

\section{Introduction}


Consider a family of probability models $P(\Y \mid \theta)$  indexed by parameter $\theta \in \Theta$ for observations $\y$, and a prior distribution $\pi$ on the parameter space $\Theta$. In a classical Bayesian approach, the posterior distribution $p(\theta \mid \Y)$ is used for statistical inference \cite{van2000asymptotic}. 

Denote the true distribution of observations $P_0$, and we consider the case $P_0 \notin \{P( \cdot \mid \theta), \, \theta \in \Theta\}$, i.e. the model is misspecified. Such case arises in many applications, particularly in complex models where the numerical evaluation of posterior distribution under the ideal probability model takes a long time to run, leading to increased use of approximate models with faster computing time. A typical example is approximating complex dependence structure by pairwise dependence only  \cite{Besag}, \cite{Varin2011}.

For well-specified regular models, the classical Bernstein - von Mises theorem 
 states that for $n$ independent identically distributed (iid) observations, for large $n$, the posterior distribution behaves approximately as a normal distribution centered on the true value of the parameter with a random shift, with both the posterior variance and the variance of the random shift being asymptotically equal to the inverse Fisher information, thus making the Bayesian inference asymptotically consistent and efficient. This was extended to locally asymptotically normal (LAN) models for $\Theta \in \mathbb{R}^p$ with fixed $p$ \cite{van2000asymptotic}, for LAN models with growing $p$, etc.

A version of Bernstein - von Mises theorem is also available for nonregular models such as locally asymptotically exponential models (\cite{IH86}, \cite{ChernoHong}), models with parameter on the boundary of the parameter space \cite{BochkinaGreen} which hold under misspecified models. Here, however, we focus on ``regular'' models where estimators are asymptotically Gaussian.

However, under model misspecification, Bayesian model is  no longer optimal, as the posterior variance does not match the minimal lower bound on the variance of unbiased estimators \cite{kleijn2012bernstein} and \cite{panov2015}). 

Therefore, the standard way of constructing a posterior distribution following the Bayes theorem may not be appropriate for particular purposes, e.g. inference or prediction \cite{muller2013risk}. Different ways  to construct a distribution of $\theta$ given $\Y$ that produces inference appropriate for the purpose of the analysis have been proposed.
A natural aim for such a method is to behave like a standard Bayesian method when the model is well specified, i.e. when $P_0 \in \{P(\cdot | \theta), \, \theta \in \Theta \}$, and to provide at least asymptotically optimal inference, from the frequentist perspective, under model misspecification.

In this review we focus on regular misspecified models.
The review is organised as follows. We start with the summary of frequentist results for regular misspecified models (Section~\ref{sec:FreqMisspec}).
 In Section~\ref{sec:ClassicalBvM}, we formulate  classical Bernstein - von Mises theorem  and in Section~\ref{sec:BvMmisspec} we discuss the analogue of Bernstein - von Mises theorem under model misspecification, particularly the conditions when this local Gaussian approximation holds. In Section~\ref{sec:GeneralisedPosterior}, we review the proposed methods to construct distribution $p(\theta \mid Y)$ that results in improved inference under model misspecification compared to the standard Bayesian approach. We conclude with discussion and open questions.

{\it Definitions.} For a vector $x \in \mathbb{R}^p$, $||x||$ denotes the Euclidean norm of $x$, and for a matrix $A \in \mathbb{R}^{p \times p}$, $||A||$ denotes the spectral (operator) norm of $A$.

\section{Frequentist results for misspecified models}\label{sec:FreqMisspec}

\subsection{Probability model}

A probability-based set up is as follows. Consider a measurable space $({\mathcal Y}, {\mathcal A})$ and let ${\mathcal P}$ be a set of probability distributions on $({\mathcal Y}, {\mathcal A})$, and assume that
 $\Y \sim P( \cdot \mid \theta)$, $\theta \in \Theta \subseteq \mathbb{R}^p$ with finite $p$ (which may or may not be allowed to grow with $n$), $\{P(\cdot \mid \theta), \theta \in \Theta\} \subset {\mathcal P}$,  where $\Y$ is $n$-dimensional random variable, $P(\y \mid \theta)$ is the probability density   function (with respect to the Lebesque or counting measure). The true distribution of $\Y$ is denoted as $P_0$, with density $p_0$.

\subsection{Best parameter}\label{sec:BestParameter}

Given the parametric family and the true distribution of the data, the ``best'' parameter is defined by
\begin{equation}\label{eq:BestParameter}
\theta^\star = \arg \min_{\theta \in \Theta} KL(  P(\cdot | \theta), P_0)
\end{equation}
where $KL(P_0 , P_1) = \int \log\left( \frac{d P_0}{ d P_1}\right) dP_0$ is the Kullback-Leibler divergence between probability measures $P_0$ and $P_1$.

Usually, it is either assumed that the model misspecification is such that $\theta^\star$ is the parameter of interest (e.g. in machine learning or quasi-likelihood approaches, this is often done by construction). However, there are alternative approaches when $\theta^\star$ differs from the parameter of interest $\theta_0$ (e.g. \cite{MillerDunson}).


\subsection{Regular models}\label{sec:RegularModels}

We consider a regular setting, under the following assumptions.
\begin{enumerate}
\item $\Theta$ is an open set.
\item Maximum of $\E \log p(\theta \mid \Y)$ over $\theta \in \Theta$ is attained at a single point $\theta^\star$ (\eqref{eq:BestParameter}).
\item $D(\theta)$ and $V(\theta)$ are  finite positive definite matrices for all $\theta$ in some neighbourhood of $\theta^\star$, where
\end{enumerate}
\begin{eqnarray}\label{eq:defDV}
  V(\theta) &=&  \E [ \nabla \log p(\theta \mid \Y)  (\nabla \log p(\theta \mid \Y))^T],\\
  D(\theta) &=& -\E  \nabla^2 \log p(\theta \mid \Y). \notag
\end{eqnarray}
Here (and throughout the paper) the expectation is taken with respect to the true distribution of the data, $\Y \sim P_0$ (which is sometimes emphasised by writing $\E_{P_0}$), and $\nabla$ is the differentiation operator with respect to $\theta$. If the model is correctly specified, $D(\theta^\star) = V(\theta^\star)$.


Typically, two main types of  models are considered in theory: independent identically distributed (iid) models: $\Y = (Y_1,\ldots,Y_n)$ with independent $Y_i$ with the same pdf or pmf $p(\cdot \mid \theta)$, and more generally locally asymptotically normal (LAN) models \cite{van2000asymptotic}.
 We give the definition of LAN models stated in \cite{kleijn2006} that applies to a misspecified model.
\begin{definition}\label{def:LAN}
 Stochastic local asymptotic normality (LAN)
condition:   given an interior   point $\theta^\star \in \Theta$ and a
rate $\delta_n\to 0$, there exist random vectors $\Delta_{n,\theta^\star}$ and a nonsingular  matrix $\D0$
such that the sequence $\Delta_{n,\theta^\star}$ is bounded in probability, and for every compact
set $K \in \mathbb{R}^p$,
$$
\sup_{ h\in K} \left| \frac{\log p_{\theta^\star +\delta_n h}(Y_1,\ldots, Y_n)}{\log p_{\theta^\star }(Y_1,\ldots, Y_n)} - h^T \D0 \Delta_{n,\theta^\star} - 0.5 h^T \D0 h \right| \to 0
$$
 as   $n\to \infty$ in (outer) $P_0^{(n)}$ -probability.
\end{definition}
For iid model and iid true distribution, with possibly misspecified density, $\delta_n = 1/\sqrt{n}$ and $\D0$ is the limit of $\delta_n^{2} D(\theta^\star)$ as $n\to \infty$  where $D(\theta^\star)$ is defined by \eqref{eq:defDV}. Also,
$$
\Delta_{n,\theta^\star} = n^{-1/2} \D0^{-1} \sum_{i=1}^n \nabla \log p (Y_i \mid \theta^\star)
$$
which has mean 0 and variance $\D0^{-1} \V0 \D0^{-1}$, called the sandwich covariance, where $\V0$ is the limit of $\delta_n^{2} V(\theta^\star)$ as $n\to \infty$, and $V(\theta^\star)$ is defined by \eqref{eq:defDV}.


\subsection{Nonasymptotic LAN condition}\label{sec:LANnonAsymp}

\cite{spokoiny2012} provides a non-asymptotic version of LAN expansion under model misspecification, and non-asymptotic bounds on consistency of the MLE and coverage of MLE-based and likelihood-based confidence regions. These conditions have been updated in \cite{panov2015}.

Define the stochastic term
$$
\zeta(\theta)=\log p(\Y \mid \theta)- \E \log p(\Y \mid \theta).
$$

\begin{enumerate}

\item There exist  $g>0$ and positive-definite $p\times p$ matrix $\Vstar$ such that for any $|\lambda|\leq g$,
$$
\sup_{\gamma \in \mathbb{R}^p} \E_{P_0}  \exp\left\{\lambda \frac{\gamma^T\nabla \zeta(\theta^\star)}{||\Vstar \gamma||}\right\} \leq e^{ \lambda^2/2},
$$
 If such matrix $\Vstar$ exists, it satisfies  $\Var(\nabla \zeta(\theta^\star)) \leq \Vstar$. \cite{panov2020} show that this holds as long as the following condition holds for some $\tilde g>0$ and $C \in (0,\infty)$:
$$
\sup_{u \in \mathbb{R}^p: \, ||u||\leq \tilde g} \E  \exp\left\{ u^T \nabla \zeta(\theta^\star) \right\} \leq C.
$$

\item There exists $\omega>0$ such that for any $||{\Dstar}^{ 1/2}(\theta-\theta^\star)|| \leq r $ and $|\lambda|\leq g$,
$$
\sup_{\gamma \in \mathbb{R}^p} \E  \exp\left\{\lambda \frac{\gamma^T(\nabla \zeta(\theta)-\nabla \zeta(\theta^\star)) }{\omega||\Vstar \gamma||}\right\} \leq e^{ \lambda^2/2}.
$$
In later work, \cite{spokoiny2020} assumes that $\nabla\zeta(\theta)$ is independent of $\theta$ by introducing an augmented model in the context of an inverse problem (the approach the author refers to as calming) thus making this condition unnecessary.

\item Conditions on $\E \log p(\Y \mid \theta)$:    $\nabla \E \log p(\Y \mid \theta^\star)=0$ and that the second derivative is continuous in the neighbourhood   $||{\Dstar}^{ 1/2}(\theta-\theta^\star)|| \leq r $:
$$
||  {\Dstar}^{-1}D(\theta)- I_p|| \leq \delta(r)
$$
   where $I_p$ is $p\times p$ identity matrix, $D(\theta )=-\nabla^2 \E \log p(\Y \mid \theta)$ and $ \Dstar =D(\theta^\star)$.  This condition is taken from \cite{panov2015}, in \cite{spokoiny2012} this condition was written in terms of  $ \E \log p(\Y \mid \theta)  -\E \log p(\Y \mid \theta^\star)$ rather than its second derivative $D(\theta)$ which are similar due to $\nabla \E \log p(\Y \mid \theta^\star)=0$.
 In later papers this condition is rewritten in terms of the moments of the third and fourth derivatives of $\E \log p(\Y \mid \theta)$ \cite{panov2020}.
\end{enumerate}

Under the conditions for the stochastic terms, using results of \cite{panov2020}, for any 
 $p\times p$ matrix $B$ such that $\trace(  B \Vstar  B^T) <\infty$, the random term can be bounded nonasymptotically as follows:
\begin{eqnarray*}
P\left(|| B (\nabla \zeta(\theta)-\nabla \zeta(\theta^\star))||\geq c_0  \omega z(B,x)\right)  \leq 2 e^{-x}
\end{eqnarray*}
where $c_0$ is an absolute constant and
$$
z(B,x) = \sqrt{\trace(B \Vstar  B^T ) + 2 x^{1/2} \trace((B \Vstar  B^T)^2) + 2x ||B \Vstar B^T||}.
$$

Note that if these conditions hold, then for $||{\Dstar}^{ 1/2}(\theta- \theta^\star)||\leq r$, with probability at least $1-2e^{-x}$,
\begin{eqnarray*}
&&\left| \log\frac{ p (\Y \mid \theta)}{ p (\Y \mid \theta^\star)} - (\theta-\theta^\star)^T   \nabla \zeta(\theta^\star)  - 0.5 (\theta-\theta^\star)^T D (\theta-\theta^\star) \right| \\
&& \quad \leq 0.5 \delta(r) r^2  + c_0 r \omega z\left({\Dstar}^{ -1/2},x\right),
\end{eqnarray*}
which is a non-asymptotic analogue of  LAN expansion under a possibly misspecified model  with $\delta_n = ||\Dstar||^{-1/2}$, $\D0 = \delta_n^2 \Dstar$ (or its limit as $\delta_n \to 0$), $h = (\theta - \theta^\star)/\delta_n$ with $||\D0^{1/2} h|| \leq \delta_n r$ and
 $\Delta_{n,\theta^\star} = \delta_n \D0^{-1}\nabla \zeta(\theta^\star)$.


\begin{example}
Consider a model with $p$-dimensional $\theta$ and   iid observations $Y_1,\ldots, Y_n$. Assume that the true observations are also iid but they may have a different true distribution  with finite positive definite $\V0 = \E \left[[\nabla \log p(Y_i\mid \theta^\star)]^T \nabla \log p(Y_i\mid \theta^\star)\right]$ and $\D0 = -\E \nabla^2 \log p(Y_i\mid \theta^\star)$. Therefore,  $\Dstar = D(\theta^\star) = n \D0$ and  $\Vstar=V(\theta^\star) = n \V0$. Denote $r_0=r/\sqrt{n}$ the radius of the local neighbourhood $||\D0^{ 1/2}(\theta- \theta^\star)||\leq r_0$.

If $\nabla \zeta(\theta)$ does not depend on $\theta$ (e.g. for $p(\cdot \mid \theta)$ from an exponential family with natural parameter), then $\omega =0$ in Condition 2 and the upper bound in the LAN condition is $0.5n r_0^3$ which tends to 0 if $r_0=o(n^{-1/3})$, or equivalently if $r=o(n^{1/6})$.

\end{example}

\subsection{Optimal variance for unbiased estimators}\label{sec:PostVar}

For regular models discussed in Section~\ref{sec:RegularModels}, the lower information bound for the variance of unbiased estimators of $\theta^\star$ when the true model is unknown,  is the sandwich covariance \cite{White1982}, \cite{Diong2017}:
$$
\Var(\hat\theta) \geq  {\Dstar}^{-1} {\Vstar} {\Dstar}^{-1}
$$
where ${\Vstar} = V(\theta^\star)$ and ${\Dstar} = D(\theta^\star)$ as defined by \eqref{eq:defDV}. This is an analogue of the Cramer-Rao inequality for regular misspecified models. In frequentist inference, the MLE is asymptotically unbiased and its variance is approximately the sandwich covariance, i.e. inference based on the MLE for misspecified models is asymptotically efficient \cite{White1982}.  

When the true model is known but a misspecified model is used, e.g. for computational convenience, it is possible in principle to achieve the smallest variance, inverse Fisher information, using a misspecified model with additional adjustment \cite{Diong2017}.  We illustrate it on the model used in  \cite{StoehrFriel} in Section~\ref{sec:Curv}.


\section{Bernstein- von Mises theorem for correctly specified models}\label{sec:ClassicalBvM}

For a correctly specified parametric model $\{P(\cdot \mid \theta), \, \theta \in \Theta \}$ with a density $p(\cdot \mid \theta)$ and a prior distribution with density $\pi(\theta)$, Bayesian inference is conducted using the posterior distribution
$$
p(\theta \mid \y) = \frac{p(\y \mid \theta) \pi(\theta)}{\int_{\Theta} p(\y \mid \theta) \pi(\theta) d\theta}, \quad \theta \in \Theta.
$$
We formulate the Bernstein - von Mises theorem in a regular setting defined in Section~\ref{sec:RegularModels}, under the additional assumption
 that  prior density $\pi(\theta)$ is continuous for in a neighbourhood of $\theta^\star$, following \cite{van2000asymptotic}.




\begin{theorem}[Bernstein - von Mises theorem]
For a well-specified regular parametric model $\{p(\y \mid \theta), \, \theta \in \Theta \}$ with $P_0 = P_{\theta^\star}$ under the regularity assumptions listed in Section~\ref{sec:RegularModels}, LAN condition with $\D0$,  for a prior density $\pi(\theta)$ continuous for in a neighbourhood of $\theta^\star$, then
$$
\sup_{ A}| P( \delta_n^{-1}(\theta - \theta^\star) \in A \mid Y_1,\ldots, Y_n) - N(A; \Delta_{n,\theta^\star},  \D0^{-1}) |  \stackrel{P_0^{\infty}}{\to} 0
$$
as $n\to \infty$, where the supremum is taken over measurable sets $A$, and $\Delta_{n,\theta^\star}$ weakly converges to $N(0, \D0^{-1})$.
\end{theorem}
 Matching variances of the posterior distribution and of the random shift $\Delta_{n,\theta^\star}$, that are equal to the inverse Fisher information, make Bayesian inference efficient asymptotically, from the frequentist perspective. As the random variable $\Delta_{n,\theta^\star}$ is bounded with high probability and $\delta_n \to 0$, this theorem also implies consistency of the posterior distribution of $\theta$.

\section{Bernstein - von Mises theorem and model misspecification}\label{sec:BvMmisspec}

\subsection{Bayesian inference under model misspecification}\label{sec:BayesMisspec}

 Given a prior distribution with density $\pi(\theta)$, the posterior distribution is constructed as the conditional distribution of the parameter $\theta$ given data $\y$ using Bayes theorem.
 A more general approach, often referred to as a Gibbs posterior distribution, is where the posterior distribution is defined using a loss function $\ell(\theta, \y)$ and a prior distribution with density $\pi(\theta)$:
$$
p(\theta \mid \y) = \frac{e^{-\ell(\theta, \y)}\pi(\theta)}{\int_{\theta \in \Theta}  e^{-\ell(\theta, \y)}\pi(\theta) d\theta}, \quad \theta \in \Theta.
$$
If the loss function $\ell(\theta, \y)$ is chosen to be $-\log p(\y \mid \theta)$, this approach leads to the usual posterior distribution. As well as differing in the interpretation, the key technical difference to the classical Bayesian approach is that function $e^{-\ell(\theta, \y)}$ does not integrate to 1 over $\y$.  This approach is used in applications where only moment conditions are available (\cite{Chib2018}, Huber function can be used as a loss for robust inference, etc.

\cite{bissiri2016general} provide a decision-theoretical justification of this approach, by showing that this distribution minimises the following loss function with respect to probability measure $\nu$ on $\Theta$,
 \begin{equation}\label{eq:GibbsPostOpt}
\int_{\Theta} \ell(\y, \theta) \nu(d\theta) + KL(\nu,\pi).
 \end{equation}
The authors argue that for iid observations, this is a Bayesian equivalent of
$$
\hat \theta = \arg \min_{\theta \in \Theta} \frac 1 n \sum_{i=1}^n \ell(y_i, \theta).
$$
In the latter approach the interest is in a point estimator whereas in the former approach the interest is in a distribution over $\theta$ given $\y$.

When applying Bayesian approach under model misspecification, the key question is whether Bayesian inference remains asymptotically efficient, i.e. whether the Bernstein - von Mises theorem holds with the posterior variance being close asymptotically to the sandwich covariance.

\subsection{Concentration}

A necessary condition for a Bernstein - von Mises - type result is to prove that the posterior distribution concentrates in the limit at the point mass at $\theta^\star$.

Consistency of the posterior distribution can be defined as follows. Given a distance $d$ between the class of probability models $\{P_\theta, \theta\in \Theta\}$ and the true distribution $P_0$, as for any $\varepsilon >0$,
$$
P( d(P_\theta, P_0) > \varepsilon \mid Y_1, \ldots, Y_n) \stackrel{P_0^{\infty}}{\to} 0 \quad \text{ as } n\to \infty.
$$
Here $P_\theta$ is the probability distribution associated with probability density (mass) function $p_\theta$.

For a misspecified model, the distance between $\{P_\theta, \theta\in \Theta\}$ and $P_0$ may be positive, so it is not always   possible to achieve consistency. However, it may be possible to prove concentration at the probability model with the best parameter $\theta^\star$:
 as for any $\varepsilon >0$,
$$
P( d(P_\theta, P_{\theta^\star}) > \varepsilon \mid Y_1, \ldots, Y_n) \stackrel{P_0^{\infty}}{\to} 0 \quad \text{ as } n\to \infty.
$$
This is referred to as posterior concentration. It is also often of interest to prove that the posterior distribution contracts at some rate (usually the corresponding minimax rate), namely that there exists a sequence $\varepsilon_n$ such that for any sequence $M_n$ growing to infinity,
\begin{eqnarray}\label{def:posteriorconsistency}
P( d(P_\theta, P_{\theta^\star}) > M_n \varepsilon_n \mid Y_1, \ldots, Y_n) \stackrel{P_0^{\infty}}{\to} 0 \quad \text{ as } n\to \infty.
\end{eqnarray}


The main paper on posterior contraction rate under model misspecification is \cite{kleijn2006}. Their results apply to nonparametric models with $\theta = P$. One of their conditions is formulated in terms of  the covering number for testing under misspecification.

\begin{definition}
Given $\varepsilon > 0$,  define
$N_t (\varepsilon, {\mathcal P}, d, P_0, P_{\theta^\star})$,  the covering number for testing under misspecification,  as the minimal number $N$ of convex sets $B_1, \ldots, B_N$ of probability
measures on $({\mathcal Y}, {\mathcal A})$ needed to cover the set $\{P \in {\mathcal P}: \, \varepsilon < d(p_\theta, p_{\theta^\star}) < 2\varepsilon\}$
such that, for every $i$,
$$
\inf_{
P\in B_i}
\sup_{0<\eta<1} -\log \E_{P_0}[p/p_{\theta^\star}]^\eta
\geq \varepsilon^2/4
$$
If there is no finite covering of this type,   the covering number is defined to be
infinite.
\end{definition}

Then, their main result (Theorem 2.1) is a typical statement on Bayesian (nonparametric) rates of posterior contraction.
Their Corollary 2.1 simplifies the statement for consistency, without the rate.
\begin{corollary}[Corollary 2.1 in \cite{kleijn2006}]
For a given model ${\mathcal P}$, prior $\Pi$ on ${\mathcal P}$  some $P_{\theta^\star}
\in {\mathcal P}$ and a semi-metric $d$ on ${\mathcal P}^2$,
assume that
\begin{enumerate}
\item $KL(p(\cdot|\theta^\star), p_{0})<\infty$
\item $\E (p(\Y|\theta)/p(\Y|\theta^\star)) <\infty$ for all $\theta \in \Theta$
\item $p(y \mid \theta^\star) >0$ for all $y\in {\mathcal Y}$
\end{enumerate}
and for  every $\varepsilon >0$,
$$
\Pi( B(\varepsilon,P_{\theta^\star}, P_0)) >0
$$
where
$$
B(\varepsilon,P_{\theta^\star}, P_0)=\left\{\theta: \, -\E_{P_0} \log \left[\frac{p_\theta}{ p_{\theta^\star}}\right] \leq \varepsilon^2, \, -\E_{P_0} \left[\frac{p_\theta}{ p_{\theta^\star}}\right]^2 \leq \varepsilon^2\right\},
$$
 and
$$
\sup_{\eta > \varepsilon} N_t (\eta, {\mathcal P}, d, P_0, P_{\theta^\star}) <\infty.
$$
Then, for every $\varepsilon >0$, as $n\to \infty$,
$$
P( d(P_\theta, P_{\theta^\star}) > \varepsilon \mid Y_1, \ldots, Y_n) \to 0 \quad \text{ as } n\to \infty.
$$
\end{corollary}

The authors also consider the case when the best approximation $\theta^\star$ is not unique. Their results are illustrated on consistency of density estimation using mixture models, and on nonparametric regression models using a convex set of prior models for the regression function.

\cite{grunwald2012safe} demonstrates that the convexity of ${\mathcal{P}}$ is crucial. \cite{GrunwaldVanOmmen} show via simulations that if $\{P_\theta, \, \theta\in\Theta\}$ is not convex, the posterior distribution does not concentrate on $\theta^\star$ but instead it concentrates on the best approximation of $P_0$ in the convex hull of the class of the parametric models $ Conv(\{P_\theta, \, \theta\in\Theta\})$:
\begin{equation}
\tilde P = \arg \min_{P \in Conv(\{P_\theta, \, \theta\in\Theta\})} KL(P_0,  P).
\end{equation}
We discuss their example in more detail in Section~\ref{sec:ExampleGrVanOmmen}.  In particular, the authors say that it is possible to achieve consistency, i.e. for the posterior distribution to converge to the point mass at $\theta^\star$ if for any $\eta \in (0,1]$,
\begin{equation}\label{eq:CondIntegralFracBayes}
\E_{P_0} \left(\frac{dP(\Y \mid \theta)}{dP(\Y \mid \theta^\star)} \right)^{\eta} \leq 1 \quad \text{ for all } \theta \in \Theta.
\end{equation}
Note that this condition is reminiscent of one of the conditions of \cite{kleijn2006} who assume that the above condition holds for $\eta=1$, and a similar condition is present in the definition of the covering numbers under model misspecification.  Further, \cite{GrunwaldMehta} relax this condition to   be upper bounded by $1+u$ for some small $u>0$. 

\cite{bhattacharya2019} study the posterior contraction rate for a particular type of semi-metric $d$ that is matched to the considered misspecified model. For such a matched semi-metric, they show that the posterior contraction rate is determined only by the prior mass condition of the posterior contraction theorem of Ghoshal and van der Vaart (2007), and does not involve the entropy condition. In their setting, the pseudo-likelihood is a power of a probability density, so  condition \eqref{eq:CondIntegralFracBayes} holds. The authors consider only examples of misspecified models with a convex parameter space.
 \cite{SyringMartinRateGibbs} study  the concentration rate  of Gibbs posteriors in a (semi-metric) $d$ under more general losses and semi-metrics, focusing on iid models and iid true distribution, and discuss a setting where $Y_i$'s are independent but not necessarily identically distributed. Their assumptions also include the condition on the prior mass of KL neighbourhood but not the entropy; they use a different additional assumption instead. See Section~\ref{sec:ConcRateGibbs} for details.


\subsection{Bernstein - von Mises - type results under model misspecification}

The first   Bernstein - von Mises - type result  under model misspecification was formulated by \cite{kleijn2012bernstein}. The authors state that for misspecified LAN models (see Definition~\ref{def:LAN}), under assumptions of Theorem 2.1 in  \cite{kleijn2006} with rate $\delta_n$,
$$
\sup_{A}| P( \delta_n^{-1}(\theta - \theta^\star) \in A \mid Y_1,\ldots, Y_n) - N(A; \Delta_{n,\theta^\star}, \D0^{-1}) |  \stackrel{P_0^{\infty}}{\to} 0
$$
  as $n\to \infty$, i.e. the posterior distribution converges to the Gaussian distribution in the total variation distance. 

\cite{panov2015} state the BvM for semi-parametric possibly misspecified models, with flat or a Gaussian prior distribution, in a non-asymptotic setting. Here we only state conditions and statement for a parametric model.

In addition to the non-asymptotic LAN assumptions stated in Section~\ref{sec:RegularModels}, the following assumptions are made.
\begin{enumerate}
\item {\it Small bias condition}:  the norm of the bias of the penalised estimator
  $ ||\theta^\star_\pi- \theta^\star||$ is small where $\theta^\star_\pi$ is defined by
\begin{eqnarray}\label{eq:BestParameterPrior}
\theta^\star_\pi &=& \arg \min_{\theta \in \Theta} [ KL(  p(\cdot | \theta) , p_0) - \log \pi(\theta)] \\
&=&  \arg \min_{\theta \in \Theta}  [- \E \log p(\Y \mid \theta) - \log \pi(\theta)].\notag
\end{eqnarray}

\item {\it Identifiability}: $||{\Dstar}^{-1} \Vstar ||\leq a^2  \in (0, \infty)$.

\item {\it Global deterministic condition}: for $||{\Dstar}^{1/2}(\theta-\theta^\star)||>r$,
$$
\E \log p(\Y \mid \theta) -\E \log p(\Y \mid \theta^\star) \geq - ||{\Dstar}^{1/2}(\theta-\theta^\star)|| b(r)
$$
with $b(r)$ growing to infinity as $r$ grows to infinity. \cite{panov2015} have a stronger condition however, it is possible to show that this condition is sufficient to bound the tail of the posterior distribution on $||{\Dstar}^{1/2}(\theta-\theta^\star)||>r$ for large $r$.

 \item {\it Global stochastic condition}:
for any $r>0$  there exists  $g(r)>0$  such that for any $|\lambda|\leq g(r)$,
 $$
\sup_{||{\Dstar}^{1/2}(\theta-\theta^\star)||\leq r } \sup_{\gamma \in \mathbb{R}^p} \E  \exp\left\{\lambda \frac{\gamma^T\nabla \zeta(\theta)}{||\Vstar \gamma||}\right\} \leq e^{ \lambda^2/2},
 $$
 where the expectation is taken with respect to $Y \sim P_0$.

\end{enumerate}



Under these assumptions, \cite{panov2015} formulated a non-asymptotic version of Bernstein - von Mises theorem for misspecified models.
\begin{theorem}[Theorem 1 in \cite{panov2015}]
Suppose that the assumptions stated in this section hold, and consider  a flat prior $\pi(\theta)=1$ for all $\theta \in \Theta$.

Then, for any measurable $A$,   with probability at least $1-4e^{-x}$,
\begin{equation}\label{eq:BvMnonasympt}
|P(\tilde D^{ 1/2}(\theta - \theta^\star -\Delta_{\star,n}) \in A \mid \Y ) - N(A;  0, I_p)| \leq  e^{\tilde\Delta(r)}-1,
\end{equation}
where $\hat\theta$ is the MLE of $\theta$, $\tilde D = D(\hat\theta)$, $\Delta_{\star,n}= \tilde D^{-1} \nabla\zeta(\theta^\star)$  and   $\tilde\Delta(r) =   r^2 (\delta(r) + 6 \omega z(\tilde D^{-1/2},x)) + 8e^{-x}$.
\end{theorem}
The authors also prove a similar result with posterior mean and posterior precision matrix instead of $\theta^\star +\Delta_{\star,n}$ and $\tilde D$.

In their Theorem 2, the authors  extend this result to a Gaussian prior $\theta \sim N(0, G^{-2})$ such that
\begin{gather*}
||G^2 {\Dstar}^{ -1}||\leq \epsilon<1/2, \quad trace [(G^2 {\Dstar}^{  -1})^2]\leq \delta^2, \\ || ({\Dstar}  +G^2)^{-1} G^2\theta^\star||\leq \beta.
\end{gather*}
Then, the posterior is approximated by the Gaussian distribution, namely, equation (\ref{eq:BvMnonasympt}) holds  with the upper bound replaced by $e^{2\Delta(r)+8e^{-x}}(1+\tau)+e^{-x}-1$ where
$$
\tau = 0.5 \left[ (1+\epsilon)(3 \beta + \epsilon z({\Dstar}^{  -1/2},x))^2 +\delta^2 \right]^{1/2}.
$$
In particular, the authors show that for a high dimensional parameter, the upper bound is small is $x$ is large and $p^3/n$ is small.
 \cite{panov2020} have shown a similar result under assumption that the stochastic term is a constant, with posterior distribution centered either at the MLE $\hat\theta$ or at the posterior mean. The authors also apply these results to nonparametric problems.

\cite{Chib2018} study Bayesian exponentially-tilted empirical likelihood posterior
distributions, which are defined by moment conditions rather than by a likelihood or
loss function.  The authors show the BvM result for well-specified and misspecified models under fairly general conditions.

 \subsection{Example: misspecified linear model}\label{sec:ExampleGrVanOmmen}

Now we consider the example of a misspecified model given in \cite{GrunwaldVanOmmen} where the Bayesian approach considered by the authors fails, and we apply theory of \cite{panov2015} to analyse it. In particular, we check whether the small bias condition holds, i.e. whether $\theta^\star_\pi$ defined by \eqref{eq:BestParameterPrior} is close to $\theta^\star$ defined by \eqref{eq:BestParameter}.

The authors considered the following linear model
\begin{eqnarray}
Y_i = \beta_0 + \sum_{j=1}^p \beta_{j} X_{ij} + \epsilon_i, \quad \epsilon_i \sim N(0,\sigma^2), \,\, i=1,\ldots,n
\end{eqnarray}
independently, with a conjugate prior distribution on its parameters:
\begin{eqnarray}
\beta = (\beta_0, \beta_{1},\ldots ,\beta_p)^T \sim N(0, c^{-1} \sigma^2 G),  \tau:=\sigma^{-2} \sim \Gamma(a,b)
\end{eqnarray}
independently, with $G=I_{p+1}$. The values of the hyperparameters were chosen to be $c=1$, $a=1$, $b=40$.

The true distribution of the data, i.e. the data generating mechanism, is as follows:
\begin{gather}
X_{ij} \sim N(0,1), \,\,  \,\,  Z_i \sim Bern(0.5),\notag\\
Y_{ij} \mid Z_i \sim  N\left(\sum_{j=1}^p \beta_{true, j} X_{ij}, \,\, \sigma_{true}^2 (1+Z_i)\right),
\end{gather}
independently for $i=1,\ldots, n$ and $j=1,\ldots,p$. The true values were taken as $n=100$, $p=40$, $\sigma_{true}^2=1/40$, $\beta_{true, j} = 0.1$ for $j=1,2,3,4$ and $\beta_{true, j} = 0$ otherwise. Note that $\beta_{true, 0}=0$, i.e. there is no intercept, and that the variance of $Y_i$ given $X$ is $\sigma_{true}^2(0.5\cdot 1 + 0.5 \cdot 2) =  1.5 \sigma_{true}^2$.

Now we work out the best parameter for this model defined by \eqref{eq:BestParameter} and the point at which posterior distribution concentrates asymptotically \eqref{eq:BestParameterPrior}, and whether they are close or not.

The log likelihood for the considered model is
$$
L(\beta, \tau) = -0.5 \tau ( Y - X \beta)^T ( Y - X \beta)  + 0.5 n \log \tau,
$$
and logarithm of the posterior distribution of $\theta=(\beta, \tau)$ is
$$
L_\pi(\beta, \tau) = L(\beta, \tau) - 0.5 c \tau \beta^T \beta + 0.5 p \log \tau +(a-1)\log \tau - b \tau.
$$
Negative Kullback-Leibler distance $KL(p_\theta, p_0)$ (up to an additive constant independent of unknown parameters), is
\begin{eqnarray}
\E L(\beta, \tau) 
&=&    0.5 n \log \tau-  0.5 \tau ||\beta - \beta_{true}||^2_2 - \frac 3 4 n \sigma^2_{true} \tau \label{eq:ExLogLik}
\end{eqnarray}
where the expectation is taken under the true model, using $\E (Y^T Y \mid X) = \beta_{true}^T X^T X \beta_{true} + 1.5 n \sigma^2_{true}$ and $\E X^T X = I_{p+1}$. Then, the best parameter, i.e. the parameter maximising  expression \eqref{eq:ExLogLik} is
\begin{eqnarray}
 \beta^\star = \beta_{true},\quad  \tau^{\star\, -1} = \sigma^{\star \, 2} =  \frac 3 2   \sigma^2_{true},
\end{eqnarray}
as stated in \cite{GrunwaldVanOmmen}.

Now we study the value of the parameters where the posterior concentrates which minimises
\begin{eqnarray}
\E L_\pi(\beta, \tau)
&=&   ( 0.5 (n + p) +  a-1) \log \tau
- \tau  b - \frac 3 4 n \sigma^2_{true}   \tau\\
&&  - 0.5 \tau (1/c+1)^{-1}  || \beta_{true}||^2_2 \\
&& -  0.5 (1+c) \tau  ||\beta -  \beta_{true}/(1+c)||^2_2
   \notag
\end{eqnarray}
and which are equal to
\begin{eqnarray}
\beta^\star_\pi &=&  (1+c)^{-1}\beta_{true},\\
\tau_\pi^{\star\, -1} &=& \sigma^{\star \, 2}_\pi = \frac{ 1.5   \sigma^2_{true}
   + 2b/n    +   c (1+c)^{-1} n^{-1} || \beta_{true}||^2_2}{ 1 + p/n  +  2(a-1)/n}. \notag
\end{eqnarray}
Hence, $(\beta^\star_\pi, \sigma^{\star \, 2}_\pi)$ is close to $(\beta^\star, \sigma^{\star \, 2} )$ if the following conditions hold:
\begin{enumerate}
\item $c=o(1)$
\item $b/n =o( \sigma^2_{true})$
\item $c (1+c)^{-1} n^{-1} || \beta_{true}||^2_2 =o( \sigma^2_{true})$
\item $p/n  =o(1)$
\item $(a-1)/n = o(1)$.
\end{enumerate}
The choice of the parameters given in \cite{GrunwaldVanOmmen} is the following:
\begin{enumerate}
\item $c=1$
\item $b/n = 0.4$, $\sigma^2_{true} = 0.025$
\item $c (1+c)^{-1} n^{-1} || \beta_{true}||^2_2 = 0.0002$,
\item $p/n  = 0.5$
\item $(a-1)/n = 0$
\end{enumerate}
i.e. conditions 3 and 5 hold whereas the remaining conditions do not hold. So, it is possible to tune hyperparameters so that all conditions, except condition 4, hold, e.g. by taking small $b$ and $c$ leading to weakly informative priors with large variances. Condition $p/n  =o(1)$ is due to the choice of the conjugate prior for $\beta$ with its prior variance proportional to the variance of the noise; if the prior variance of $\beta$ does not depend on the noise variance, then this condition is not necessary.

For instance, it is easy to show using the same technique, that considering a non-conjugate prior $\beta \sim N(0, \tau_0^{-1}\Sigma_0)$ and $\tau \sim \Gamma(a,b)$, with $||\Sigma_0||=1$, under the following conditions
\begin{enumerate}
\item $\tau_0=o(1)$
\item $b/n =o( \sigma^2_{true})$
\item $(a-1)/n = o(1)$
\end{enumerate}
leads to $\theta_\pi^\star$ being close to $\theta^\star$. These conditions are satisfied e.g. with small $\tau_0$, small $b$ and $a=1$, provided $\sigma^2_{true}$ is not much smaller than $1/n$.

\section{``Optimising'' Bayesian inference under model misspecification}\label{sec:GeneralisedPosterior}

\subsection{Asymptotic risk of parameter estimation under a misspecified model}


\cite{muller2013risk}  showed that the asymptotic frequentist risk associated
with misspecified Bayesian estimators is inferior to that of an artificial posterior which is
normally distributed, centred at the maximum likelihood estimator and with the sandwich covariance
  matrix.

This provided theoretical justification for constructing a (quasi-) posterior  distribution based on a misspecified model such that its posterior variance is approximately the sandwich covariance. Several such approaches have been used that we discuss below.


\subsection{Composite likelihoods}



Composite likelihoods (also known as pseudo-likelihoods) have been studied by \cite{LindsayCompLik}, and they are defined as follows.
Denote by $\{A_1,\ldots, A_K\}$ a set of marginal or conditional events with associated likelihoods $L_k(\theta; y) \propto P(y \in A_k; \theta)$. Then, a composite likelihood is the weighted product
$$
L_c(\theta; y) = \prod_{k=1}^K [L_k(\theta; y)]^{w_k},
$$
where $w_k$ are nonnegative weights to be chosen. It is often used to simplify the model for dependence structure in time series and in spatial models, with one of the most famous examples given by \cite{Besag} of approximating spatial dependence by the product of  conditional densities of a single observation given its neighbours. Selection of unequal weights to improve efficiency  in the context of particular applications and a review of frequentist inference for this approach is discussed by \cite{Varin2011}. For the discussion of connection of the choice of weights with the Bayesian inference under empirical likelihood see \cite{SCHENNACHELBayes}. A typical example is when $L_k(\theta; y)$ is the marginal likelihood for $y_k$ (with $K=n$). 
 Unless more information is available, it is generally difficult to estimate individual weights from the sample however this formulation gave rise to a number of approaches with randomised weights $(w_1,\ldots, w_n)$.

The idea of composite likelihood is used to sample the powers (weights of the contributions of individual samples $w_i$) from some probability distribution. The typical choice of a joint Dirichlet distribution for the weights corresponds to Bayesian bootstrap and is discussed in Section~\ref{sec:BootstrapPost}. Other choices of weights and their effect on the corresponding posterior inference are discussed in \cite{PriorWeights}. As far as I am aware, currently there are no BvM results for other randomisation schemes, apart from a joint Dirichlet distribution.

Choosing the same weight $w_k = w$ for all $k$ leads to fractional or tempered posterior distributions discussed in Section~\ref{GibbsDistr}.

\subsection{Generalised (Gibbs) posterior distribution}\label{GibbsDistr}

\subsubsection{Definition and interpretation}

Let $\ell (\y, \theta)$ be some loss function.  Then, generalised posterior distribution is given by
\begin{equation}\label{def:GibbsPosterior}
p_\eta(\theta \mid y_1,\ldots y_n) = \frac{\exp\{ - \eta  \ell (\y, \theta) \}\pi(\theta)}{\int \exp\{ - \eta \ell(\y, \theta) \pi(\theta) d\theta }
\end{equation}
where $\eta$ is the parameter that adjusts for misspecification. This parameter is called the learning rate (in machine learning), inverse temperature.
Taking $\ell(\y, \theta)= \sum_{i=1}^n\ell_i(y_i, \theta)$ for some loss $\ell_i$ associated with observation $y_i$ given parameter $\theta$ corresponds to the assumption that observations $y_i$ are  independent. Taking $\ell(\y, \theta) = -\log p(\y \mid \theta)$ and $\eta=1$ leads to the classical Bayesian inference. 
 Different functions $\ell$ may be used for different types of model misspecification and different inference purposes, e.g.  Huber function for robust  parameter estimation.

 This is also known as a Gibbs posterior in Bayesian literature, exponential weighting in frequentist literature \cite{DalalyanTsybakovPACBayes}, typically with $\ell_i = ||y_i-\theta||_2^2$, and it is used as a model for PAC-Bayesian approach in machine learning. Lately it has also been referred to as a fractional posterior and as a tempered posterior.

\cite{GrunwaldVanOmmen} argue that if there exists $\bar \eta \leq 1$ such that for all $0<\eta\leq \bar\eta$
$$
\int p_0(\y)\left(\frac{p(\y \mid \theta)}{p(\y \mid \theta^\star)} \right)^{\eta} d\y \leq 1 \quad \text{ for all } \theta \in \Theta,
$$
then the generalised posterior with $\eta < \bar\eta$ is asymptotically consistent, i.e. converges to the point mass at $p_{\theta^\star}$.
 For  $\eta$ such that condition \eqref{eq:CondIntegralFracBayes} holds, the authors interpret the generalised  posterior as a posterior distribution based on the reweighted true likelihood:
$$
\tilde p_{\eta, \theta^\star}(\y \mid \theta) = p_0(\y)\left(\frac{p(\y \mid \theta)}{p(\y \mid \theta^\star)} \right)^{\eta}
$$
which is interpreted as a density on the probability space augmented by an unobserved event if this function integrates to a positive value less than 1.
This is due to the following: if this density was used as a density of $y$ given $\theta$ to construct the likelihood, then this would correspond to a correctly specified model, since  the KL distance between $\tilde p_{\eta, \theta^\star}(\y \mid \theta)$ and $p_0$ is minimised at $\theta^\star$ and $\tilde p_{\eta, \theta^\star}(\cdot \mid \theta^\star) = p_0(\cdot)$, and the corresponding posterior would be a proper posterior and it coincides with the generalised posterior. This is done for interpretation only, as it is not possible to use $\tilde p_{\eta, \theta^\star}$ for inference in practice due to unknown $p_0$ and $\theta^\star$.

In the following section we discuss known results about concentration of the Gibbs posterior distribution.





\subsubsection{Concentration and posterior contraction rate}\label{sec:ConcRateGibbs}



\cite{bhattacharya2019} study the posterior contraction rate for a particular type of semi-metric $d$ that is matched to the considered misspecified model. They consider generalised Bayesian approach
$$
p_\eta(\theta \mid \y) = \frac{[p(\y \mid \theta)]^\eta \pi(\theta)}{ \int [p(\y \mid \theta)]^\eta \pi(\theta) d\theta}
$$
with $p(y \mid \theta)$ being a density of a probability measure with respect to some measure $\mu$,
and the semi-metric based on Renyi divergence with {\it matching} index $\eta$:
$$
D_{\eta}^{(n)} (\theta, \theta^\star) = -\frac{1}{1-\eta} \log A_{\eta}^{(n)} (\theta, \theta^\star)
$$
where $A_{\eta}^{(n)} (\theta, \theta^\star)$ is   the integral defined in \eqref{eq:CondIntegralFracBayes} for all $n$ observation $y=(y_1,\ldots, y_n)$:
$$
A_{\eta}^{(n)} (\theta, \theta^\star) = \E_{P_0}  \left[  \left(\frac{p(\Y \mid \theta)}{p(\Y \mid \theta^\star)} \right)^{\eta}\right].
$$
The authors show that since $p(y \mid \theta)$ is a density of a probability measure, condition $A_{\eta}^{(n)} (\theta, \theta^\star)\leq 1$ (equation \eqref{eq:CondIntegralFracBayes}) holds. Also, the authors show that for $\eta \to 1-$, the generalised posterior converges to the corresponding posterior  distribution. 
 Therefore, their approach is not shown to apply to so called Gibbs posteriors where other loss functions (rather than a negative log density) can be used to specify the (pseudo-)likelihood.

Following \cite{GrunwaldVanOmmen}, one may argue that in the setting considered by \cite{bhattacharya2019}, it is not necessary to use   $\eta<1$ to adjust the inference to achieve posterior consistency (it may be necessary e.g. to achieve asymptotic efficiency).
Under the assumptions of \cite{bhattacharya2019}, $A_{\eta}^{(n)} (\theta, \theta^\star) \leq 1$ for all $\theta$ and $\eta \in (0,1)$, and, due to convergence of the generalised posterior to the posterior as $\eta \to 1-$, the posterior distribution (with $\eta=1$) is asymptotically consistent.



\cite{SyringMartinRateGibbs} study  the concentration rate  of a Gibbs posterior \eqref{def:GibbsPosterior} in (semi-metric) $d$ under the following fairly general assumptions. The authors focus on iid models and iid true distribution, and discuss a setting where $Y_i$'s are independent but not necessarily identically distributed.

{\it {\bf  Condition 1.} There exist  $\bar\eta, K, r >0$  such that for all $\eta \in (0,\bar\eta)$  and for all sufficiently small $\delta > 0$, for  $\theta \in \Theta$,
$$
d(\theta, \theta^\star) > \delta \Rightarrow \log \E \exp(-\eta( \ell(Y_1; \theta)-  \ell(Y_1; \theta^\star; Y_1))) <  -K \eta \delta^r.
$$
}

{\it {\bf  KL neighbourhood condition.} For a sequence $(\varepsilon_n)$ such that $\varepsilon_n \to 0$ and $n \varepsilon_n^r \to \infty$ as $n\to \infty$, there exists $C_1 \in (0,\infty)$ such that for all $n$ large enough,
$$
\log \Pi(B_{KL}( \varepsilon^r) ) > - C_1 n\varepsilon_n^r,
$$
where
\begin{eqnarray*}
B_{KL}(R) = \left\{\theta\in \Theta: \,   -\E [\ell(\Y; \theta^\star) -  \ell(\Y; \theta)] \leq R \, \right.\\
     \quad \quad   \left.    \Var [\ell(\Y; \theta^\star) -  \ell(\Y; \theta)]\leq R \right\}.
\end{eqnarray*}
}

\begin{theorem}[Theorem 3.2, \cite{SyringMartinRateGibbs}]
Under Condition 1 and KL neighbourhood condition,  for a fixed $\eta$, the Gibbs
posterior distribution defined by \eqref{def:GibbsPosterior} satisfies   \eqref{def:posteriorconsistency} with asymptotic concentration rate $\varepsilon_n$.
\end{theorem}
The authors show that this also holds for $\eta_n \to 0$ as  long as $\eta_n n\varepsilon_n^r \to \infty$ and in the KL neighbourhood condition $C_1 n\varepsilon_n^r$ is replaced by $ C_1\eta_n n\varepsilon_n^r$. The also how that this holds for a random $\hat\eta$ as long as with high probability $c^{-1}\eta_n\leq \hat\eta \leq c\eta_n$ for $\eta_n\to 0$ and some $c\geq 1$.

The authors also discuss that condition \eqref{eq:CondIntegralFracBayes} can be relaxed to hold on $\Theta_n =\{\theta \in \Theta: \, ||\theta|| \leq \Delta_n\}$ for a sequence $(\Delta_n)$ increasing to infinity, under stronger conditions (see Theorem 4.1 in \cite{SyringMartinRateGibbs}). The authors also discuss that conditions of this theorem are related to the entropy condition of \cite{kleijn2006} and verify this condition for convex $\ell$ as a function of $\theta$.

In the iid setting, Condition 1 combines several conditions  of \cite{spokoiny2012} for a single observation $Y_1$ since
\begin{gather}\label{eq:logEexpEx3}
\log \E \exp(-\eta  (\ell(Y_1; \theta)-\ell(Y_1; \theta^\star))) = \notag\\
 -\eta(\E \ell(Y_1; \theta) - \E \ell(Y_1; \theta^\star))
  +  \log \E \exp(\eta (\zeta_1( \theta)-\zeta_1(\theta^\star))),
\end{gather}
where $\zeta_1(\theta) = \E \ell(Y_1, \theta)-\ell(Y_1, \theta)$, except that the authors assume that this condition holds for all $\theta \in \Theta$ whereas in \cite{spokoiny2012} the conditions are split into local (in a neighbourhood of $\theta^\star$) and global (for all $\theta \in \Theta$), with the global conditions being weaker.

As the authors discuss, their Condition 1 can hold if $\zeta_1(\theta)-\zeta_1(\theta^\star)$ has sub-Gaussian tails, and if for $\eta$ small enough the first term (which is negative) in \eqref{eq:logEexpEx3}  is sufficiently greater in absolute value than the second term. More specifically, assume that there exist $r, K_1 >0$ such that
$$
 -[\E \ell(Y_1; \theta) - \E \ell(Y_1; \theta^\star)] < - K_1 [d(\theta, \theta^\star)]^r, \quad \theta, \theta^\star \in \Theta,
$$
and that the sub-Gaussian tail condition  holds with some $b>0$
$$
\log \E \exp(\eta  (\zeta_1(\theta)-\zeta_1(\theta^\star))) \leq b \eta^2 ||\theta - \theta^\star||_2^2/2
$$
which can be verified through Conditions 1 and 2  \cite{spokoiny2012}, and there exist $K_2 >0$ such that $||\theta - \theta^\star||_2^2\leq K_2 [d(\theta, \theta^\star)]^r$ for all $\theta, \theta^\star \in \Theta$. Then,
\begin{eqnarray*}
\log \E \exp(-\eta (\ell(Y_1; \theta)-\ell(Y_1; \theta^\star))) &\leq&  - \eta  [d(\theta, \theta^\star)]^r [ K_1-  b  \eta  K_2 /2]\\
 &<& - K \eta  [d(\theta, \theta^\star)]^r
\end{eqnarray*}
 if $\eta <\bar\eta = \min(1+o(1), 2 K_1/(b  K_2 ))$ and $K= K_1-  b  \bar\eta  K_2 /2$. Since the inequality   $\eta < \bar\eta$ is strict, as long as $2 K_1/(b  K_2 ) >1$, we can take $\eta=1$.

The authors suggest that case $r=2$ corresponds to regular problems, i.e. where $\nabla  \E \ell(\theta^\star)=0$ and $\nabla^2  \E \ell(\theta^\star)$ is positive and continuous in the neighbourhood of $\theta^\star$, and the sub-Gaussian tails condition, whereas  nonregular problems may require other values of $r$, e.g. $r=1$ if $\theta^\star$ is on the boundary of the parameter space \cite{BochkinaGreen}, or if there is a finite jump at $\theta$ \cite{ChernoHong}. We illustrate the latter on an example.

\begin{example}\label{eq:ExCond1Jump} Now we check if Condition 1 holds for a density with jump. Consider a density $p(y| \theta)$ that is 0 for $y<\theta$ and the right hand side limit $\lim_{y \to \theta+} p(y| \theta)= \lambda >0$, for instance with $p(y| \theta) =  e^{- (y-\theta)}$ for $y>\theta$, and the true density   $p_{0}(y)$ such that $p_0(y)=0$ for  $y<\theta_0$ and $\lim_{y \to \theta_0+} p_0(y)= c_0 >0$. Then,
\begin{eqnarray*}
\E \log p(Y \mid \theta) &=& -\int_{\theta_0}^{\infty} (y-\theta) p_0(y) dy + \log (0) I(\theta > \theta_0)\\
 &=& \theta- \E Y  + \log (0) I(\theta > \theta_0)
\end{eqnarray*}
which is minimised at $\theta^\star = \theta_0$. This implies that for $\theta \in \Theta_0 = (-\infty, \theta^\star]$,
$$
\E \log p(Y \mid \theta) - \E \log p(Y \mid \theta^\star) = \theta - \theta^\star
$$
So if $d(x,y)= |x-y|$  then this holds with $K_1=1$ and $r=1$.

The stochastic term is $\zeta(\theta) = Y -\E Y$ for $\theta \in \Theta_0$, and for $\theta \in \Theta_0$
$$
\log \E \exp(\eta (\zeta(\theta)-\zeta(\theta^\star))) = 0,
$$
and Condition 1 holds for $\theta \in (-\infty, \theta^\star]$ with $K=1$ and $r=1$.
\end{example}



\subsubsection{Estimation of $\eta$}

There are various approaches to estimation of $\eta$ that lead to the posterior distribution concentrating at $\theta^\star$, e.g. \cite{GrunwaldMehta} and \cite{CooleyCompLik};  see a review \cite{ReviewEstEta}. Here I will give a very brief discussion.  
  There are two key issues: firstly, this parameter models misspecification so it cannot be estimated in a usual Bayesian way (e.g. by putting a hyperprior), and secondly, a relevant estimator depends on the aim of the inference.

  1. When predictive inference is of interest, the Safe-Bayes estimator of \cite{grunwald2012safe} further explored in \cite{GrunwaldMehta}, may be appropriate:
  $$
  \hat\eta = \arg \min \left[ - \sum_{i=1}^n\E_{} \log p(y_i\mid  \theta) p_\eta(\theta \mid y_{1:(i-1)}) d\theta \right]
  $$
where $p_\eta(\theta \mid y_{1:(i-1)})$ is the generalised posterior distribution with parameter $\eta$ based on $i-1$ samples (if $i=1$ then it is the prior distribution). 

 2. Now we discuss estimators of $\eta$ when estimation of $\theta$ is of interest, in particular frequentist coverage of  credible posterior regions  $C_{\alpha}$ such that $P_\eta(\theta \in C_{\alpha} \mid \Y) =1-\alpha$.


  Under the conditions of Gaussian approximation of the posterior, if asymptotic  coverage of asymptotic credible balls $C_{\alpha} $ is of interest, then it is sufficient to check that asymptotic credible balls
  $$
  (\theta-\hat\theta)^T {\Dstar}_\eta (\theta-\hat\theta)^T\leq \chi_p^2(\alpha)
  $$
   are inside the frequentist confidence balls with the sandwich covariance
   $$
    (\theta-\hat\theta)^T {\Dstar}  {\Vstar}^{-1}  {\Dstar}   (\theta-\hat\theta)^T\leq \chi_p^2(\alpha),
   $$
i.e. it is sufficient to check that the largest eigenvalue of $\Dstar_\eta = \eta \Dstar$ is not smaller than the largest eigenvalue of ${\Dstar}  {\Vstar}^{-1} {\Dstar}$.

Therefore, $\eta$ is chosen so that the largest eigenvalue of the posterior precision matrix $\eta {\Dstar}$ matches that largest eigenvalue of the sandwich precision matrix ${\Dstar}  {\Vstar}^{-1} {\Dstar}$, i.e. the ``oracle'' value is
$$
\eta^\star =  \frac{||{\Dstar}   {\Vstar}^{-1} {\Dstar}||}{|| {\Dstar} ||},
$$
and it can be estimated if estimates of $\Vstar=V(\theta^\star)$ and $\Dstar=D(\theta^\star)$ are available.

  \cite{holmes2017assigning} used the Fisher information number to calibrate this parameter:
  $$
\eta^\star =  \frac{\trace({\Dstar}   {\Vstar}^{-1} {\Dstar})}{\trace({\Dstar})}.
  $$
As the authors say, it is the sum of the marginal Fisher information for each dimension, which can be used as a summary for the amount of information in a sample about  parameters.

\cite{PauliVentura2012} propose $\hat\eta = \trace(V^{-1}(\hat\theta_c) D(\hat\theta_c))/\dim(\theta)$ which asymptotically is the average of the mutual eigenvalues of $D(\hat\theta_c))$ with respect to $V(\hat\theta_c)$.


\begin{remark}
Suppose that conditions of \cite{panov2015} hold for some pseudo likelihood $p(\y \mid \theta)= \prod_{i=1}^n p(y_i \mid \theta)$ and $\theta \in \Theta \subseteq \mathbb{R}$. Then, these conditions hold for $p(\y \mid \theta)^{\eta}= \prod_{i=1}^n p(y_i \mid \theta)^\eta$ with $\Vstar_\eta =  \eta^2 \Vstar$ and $\Dstar_\eta =  \eta \Dstar$,   $r_\eta = r \eta$ and $b_\eta$ such that $b_\eta(r_\eta) = b(r)$. 

Therefore, the posterior variance  ${\Dstar}^{-1}_\eta$ coincides with the sandwich variance if $\eta {\Dstar} = {\Dstar}^{ 2} {\Vstar}^{ -1}$, i.e. if $\eta = {\Dstar} {\Vstar}^{ -1}$. See also \cite{Grunwald2020} in the context of linear regression.
\end{remark}

\subsection{Nonparametric model for uncertainty in $p_0$ and bootstrap posterior}\label{sec:BootstrapPost}

\subsubsection{Nonparametric model and connection to bootstrap}

\cite{LyddonLossBootstrapBvM} proposed  to take into the account uncertainty about the parametric model by modelling the distribution of the data nonparametrically, e.g. using a Dirichlet process prior with the base model being the considered parametric model:
$$
y_i \sim F,\quad i=1,\ldots,
F \sim DP(\alpha, p(\cdot\mid \theta)).
$$

For iid observations $y_i$, when $\alpha \to 0$, this approach corresponds to Bayesian bootstrap \cite{RubinBootstrap}, with the following sampling of $(\theta^{(j)})_{j=1}^B$ from the bootstrap posterior:
\begin{gather}\label{eq:LossBootstrap}
\theta^{(j)} =  \theta(F_j) \text{ with } F_j(x) = \sum_{j=1}^n \alpha_{ji} \delta_{y_i}(x) \\
\text{and } \alpha_{j} = (\alpha_{j1},\ldots \alpha_{jn}) \sim Dirichlet(1,\ldots,1),
\end{gather}
where
\begin{equation}\label{eq:LossBootstrapTheta}
\theta(F) = \arg\min_{\theta \in \Theta} \int \ell( \theta, y) dF(y)
\end{equation}
with $\ell(\theta, y)=-\log p(y \mid \theta )$. More generally, for a possibly different loss function, the authors refer to this as the  loss-likelihood (LL) bootstrap approach. The authors argue that using this procedure induces a prior distribution on $\theta$ defined as  $P(\theta \in A)= P(F: \, \theta(F) \in A)$.

\subsubsection{Asymptotic normality of bootstrap posterior}

\cite{LyddonLossBootstrapBvM} show that the sample from the loss-likelihood bootstrap has asymptotically normal distribution with sandwich covariance matrix, weakly, under the following  assumptions.
\begin{enumerate}

\item $\Theta$ is a compact and convex
  subset of a $p$-dimensional Euclidean space.

\item The loss function $\ell: \, \Theta \times \mathbb{R} \to \mathbb{R}$ is a measurable bounded from
below function, with
$$
\int \ell( \theta, y) p_0(y) dy < \infty \quad \text{ for all } \theta \in \Theta
$$

\item (Identifiability). There exists a unique minimizing parameter value
$$\theta^\star =\arg \min_{\theta \in \Theta} \int \ell(\theta, y) p_0(y) dy,
$$
and for all $\delta>0$ where exists $\epsilon>0$ such that
    $$
\lim \inf_n  P(\sup_{|\theta-\theta^\star|>\delta} \frac 1 n \sum_{i=1}^n [\ell(\theta, y_i) - \ell(\theta^\star, y_i) ] >\epsilon ) =1
$$
\item Smoothness of loss: there exists an open ball $B$ containing $\theta^\star$ such that $$
\E[|\nabla^k_{I_k} \ell(\theta, Y)|]  <\infty   \text { and } \E[|\nabla_{j} \ell(\theta, Y) \nabla^2_{km} \ell(\theta, Y)|] <\infty
    $$
for $k=1,2,3$ and for all  corresponding indices  $I_k$, i.e. $I_1 \in \{1:p\}$, $I_2 \in \{1:p\}^2$, $(j,k,m) = I_3 \in \{1:p\}^3$, where $Y \sim p_0$.

\item For $\theta \in B$, the corresponding information matrices $V(\theta)$ and $D(\theta)$ are positive definite with all elements being finite, where
\begin{eqnarray*}
V(\theta) &=& \int \nabla \ell(\theta, y) \nabla^T \ell(\theta, y) p_0(y) dy, \\
D(\theta) &=& \int \nabla^2 \ell(\theta, y) p_0(y) dy.
\end{eqnarray*}
\end{enumerate}

\begin{theorem}[Theorem 1 in \cite{LyddonLossBootstrapBvM}]
Let $\tilde\theta_n$ be a loss-likelihood bootstrap sample of a parameter defined by \eqref{eq:LossBootstrap} and \eqref{eq:LossBootstrapTheta}  with
loss function $\ell$, given $n$ iid observations $(y_1,\ldots, y_n)$, and let $P_{LL}$ be its probability measure. Under  the above assumptions, for any Borel set $A \in \mathbb{R}^p$, as $n\to \infty$,
$$
P_{LL}(n^{1/2}(\tilde \theta_n - \hat\theta_n \in A)) \to P(Z \in A)
$$
where $Z \sim N_p(0, \Dstar {\Vstar}^{-1} \Dstar)$, $\hat\theta_n = \arg\min_{\theta \in \Theta} \frac 1 n \sum_{i=1}^n \ell(\theta, y_i)$ and
$$
\Vstar = V(\theta^\star), \quad \Dstar = D(\theta^\star, y).
$$
\end{theorem}
Therefore, the inference approach based on the loss-likelihood bootstrap is asymptotically efficient in the case the true distribution of the observations is unknown. Strictly speaking, this is not a Bernstein - von Mises theorem, since the convergence is not in the total variation distance, and hence it does not guarantee approximation of $\Pi(\theta \in A \mid \Y)$ by the corresponding Gaussian probabilities for all Borel sets $A$. Also, assumption of compactness of $\Theta$  is not present in other results on posterior concentration, so it should be possible to relax this assumption.

Another interesting problem is how to modify this approach to take into the account a given a prior $\pi$ that results in coherent and efficient inference about parameter $\theta$. \cite{NewtonWeightedBootstrap} proposed such a solution, by replacing the loss function in the optimisation problem \eqref{eq:LossBootstrapTheta} by the loss function penalised by negative log prior, however for their choice of weights, the authors give a heuristic argument that their method approximates the target posterior with posterior covariance ${\Dstar}^{-1}$ rather than with the sandwich covariance.






\subsubsection{Other bootstrap-based approaches }

Another approach is called bagged posterior or ``BayesBag'' which applies bagging proposed by\cite{BreimanBagging} to the Bayesian posterior \cite{WaddellBagging}. The idea is to select subsets of data as in bootstrap, compute posterior distribution for each of these subsets of data and average these posteriors. Formally, the bagged posterior is defined by
$$
p_{BayesBag}(\theta \mid \Y) = \frac 1 {|I|} \sum_{\Y_{(i, N)} \in I} \pi(\theta \mid \Y_{(i, N)})
$$
of the original data $\Y = (Y_1, \ldots, Y_n)$ and  bootstrap data sets $\Y_{(i, N)}$ of size $N$ as the observed data. \cite{BayesBag2020} show that under  a range of conditions, for iid true distribution of the data and iid model,  bagged posterior distribution of $\sqrt{n}(\theta - \E (\theta \mid \Y))$ converges weakly to a Gaussian distribution centered at 0 with covariance matrix ${\Dstar}^{-1}/c + {\Dstar}^{-1}\Vstar {\Dstar}^{-1}/c$ where $c = \lim (N/n)$. Hence, this approach does better than the usual posterior distribution with e.g. $N = n-n_0$ for some small finite constant $n_0$ leading to $c=1$, however it is still not efficient.

\subsection{Curvature adjustment}\label{sec:Curv}

 The generalised posterior approach uses a single parameter $\eta$ to adjust for model misspecification. In general, it is possible to use this approach to obtain variance adjustment - and hence asymptotically optimal and valid posterior inference - only for one-dimensional parameter $\theta$. In the case of higher dimensions, \cite{CooleyCompLik} proposed to use curvature adjustment in the following way. For a possibly misspecified parametric family $\{p(\dot | \theta), \, \theta \in \Theta\}$ and prior $\pi(\theta)$ ,  consider the following family of posterior distributions:
$$
\pi_A(\theta  | y  ) = p(A \theta \mid y)\propto p(y\mid A\theta)\pi(A\theta).
$$

Then, the idea is to find an estimator of the`` oracle'' matrix $A$ in this class of admissible transforms determined by the condition that the posterior variance of $A\theta$ is optimal, i.e. under the condition
$$
A^T {\Dstar}^{-1} A = {\Dstar}^{-1} \Vstar {\Dstar}^{-1}
$$
in the case the true parametric model is unknown and $\Vstar= V(\theta^\star)$ and $\Vstar= V(\theta^\star)$ defined by \eqref{eq:defDV}, and under condition $ A^T {\Dstar}^{-1} A = \Vtrue^{-1} $ if the true parametric model is known with $\Vtrue$ being the Fisher information under the true parametric model. For the case of the unknown parametric model, $\Dstar$ can be estimated by the posterior precision matrix of $\theta$, and estimation of $\Vstar$ is usually more challenging.

  \cite{StoehrFriel} applied an affine version of the transform, i.e. they considered $p(A \theta + b \mid y)$   misspecified models with known true parametric model, estimating $b$ and $A$ so that the posterior mean and the posterior variance of this distribution coincide with  the posterior mean and the posterior variance under the posterior distribution with the true parametric model.

\section{Discussion and open questions }

The approach of \cite{panov2015} and \cite{panov2020} allows to address numerically the approximation properties of misspecified Bayesian inference and to verify whether it is close to being efficient, or whether a further adjustment is needed. While the authors have the assumption of a flat or Gaussian prior, for many model it is fairly straightforward to extend this to a larger class of continuous priors, in some cases with a continuous second derivative of the log likelihood.

There are many other interesting aspects of inference under model misspecification that are not considered here, such as optimality of predictive inference,  model selection, etc.

It would be interesting to explore the connection between PAC-Bayesian inequalities and the conditional distribution $p(\theta \mid y)$ defined as the solution of the optimisation problem \eqref{eq:GibbsPostOpt}.

Other interesting approaches include BvM for Variational Bayes under model misspecification \cite{VariationalBayesBvM},
 BvM for median and quantiles under classical and Gibbs posterior \cite{BvMforMediansQuantiles}, \cite{MVquantiles}. Interestingly, \cite{zhang2021}
 show that Bayesian neural networks show inconsistency similar to that discussed in \cite{GrunwaldVanOmmen}, applying variational Bayes leads to BNN becoming consistent; it would be interesting to study whether it is possible to achieve asymptotic efficiency.
 Another version of robust Bayes-like estimation is proposed by \cite{BaraudBirge} that does not involve Kullback - Leibler distance but is based only on Hellinger distance between the true distribution and the parametric family.

Construction (asymptotically) efficient more general  Bayesian inference under model misspecification (which is also computationally tractable) is a very active research area, with several promising solutions such as bootstrap posterior and curvature adjustment, however there is still no general unifying framework to encompass these approaches or to provide a coherent general framework. Fractional posterior allows a potentially simpler  procedure for model correction which involves a single tuning parameter even if the parameter is multivariate which may result in conservative inference which can be sufficient for some problems but it is unlikely to be efficient in general. Linear curvature adjustment appear to work in practice and it is applicable to the models with no independence structure but there is no decision - theoretic justification for this is available yet; such justification is likely to involve geometry of the model space and its local linear adjustment. The open question in bootstrap-based posterior inference is the use of a given prior and its extension to data without independence structure which is likely to come from its Bayesian nonparametric interpretation.  \cite{Knoblauch2021} proposed an approach combining generalised variational inference, PAC-Bayes and other approaches into a single principled framework; they give conditions for consistency of their approach but not for efficiency. \cite{Holmes2021} propose a novel view to constructing a generalised posterior distribution, so it would be interesting to study its efficiency.


\section*{Acknowledgement.}
This review was in part motivated by the discussion of the author with Peter Gr\"unwald, Pierre Jacob and Jeffrey Miller during a Research in Groups meeting sponsored by the International Centre for Mathematical Sciences in Edinburgh, UK.


\end{document}